\documentclass{ifacconf}

\usepackage{amssymb}
\usepackage{graphicx}
\usepackage{epstopdf}
\usepackage{epsfig}
\usepackage{amsmath}
\usepackage{times}
\usepackage{verbatim}
\usepackage[dvipsnames]{xcolor}
\makeatletter
\let\old@ssect\@ssect 
\makeatother
\usepackage{natbib}
\usepackage{hyperref}
\makeatletter
\def\@ssect#1#2#3#4#5#6{%
  \NR@gettitle{#6}
  \old@ssect{#1}{#2}{#3}{#4}{#5}{#6}
}
\makeatother

\usepackage{enumitem}
\usepackage{tabularx}

\makeatletter
\def\endfigure{\end@float}
\def\endtable{\end@float}
\makeatother

\let\ifacconfcaptionwidth\captionwidth
\usepackage[caption=false]{subfig}
\let\captionwidth\ifacconfcaptionwidth

\begin{document}
\begin{frontmatter}
\title{Optimization-based Fault Mitigation for Safe Automated Driving}

\author[First]{Niels Lodder}
\author[Second,Third]{Chris van der Ploeg}
\author[First]{Laura Ferranti}
\author[Second,Third]{Emilia Silvas}%
\address[First]{Department of Cognitive Robotics, Delft University of Technology, 2628 CD Delft, The Netherlands}%
\address[Second]{TNO - Integrated Vehicle Safety, 5708 JZ Helmond, The Netherlands}
\address[Third]{Department of Mechanical Engineering, Eindhoven University of Technology, 5612 AZ Eindhoven, The Netherlands}%
\begin{abstract}
With increased developments and interest in cooperative driving and higher levels of automation (SAE level 3+), the need for safety systems that are capable to monitor system health and maintain safe operations in faulty scenarios is increasing.
A variety of faults or failures could occur, and there exists a high variety of ways to respond to such events. Once a fault or failure is detected, there is a need to classify its severity and decide on appropriate and safe mitigating actions.
To provide a solution to this mitigation challenge, in this paper a functional-safety architecture is proposed and an optimization-based mitigation algorithm is introduced. This algorithm uses nonlinear model predictive control (NMPC) to bring a vehicle, suffering from a severe fault, such as a power steering failure, to a safe-state. 
The internal model of the NMPC uses the information from the fault detection, isolation and identification to optimize the tracking performance of the controller, showcasing the need of the proposed architecture.
Given a string of ACC vehicles, our results demonstrate a variety of tactical decision-making approaches that a fault-affected vehicle could employ to manage any faults. Furthermore, we show the potential for improving the safety of the affected vehicle as well as the effect of these approaches on the duration of the manoeuvre.
\end{abstract}
\end{frontmatter}
\thispagestyle{empty}
\pagestyle{empty}

\section{INTRODUCTION}
Cooperative and automated driving (e.g., platooning) have been widely researched in the past decades, showing their effects on reducing workload and stress of the drivers \cite{Heikoop2017EffectsStudy}, but also on society. Driving in a platoon can increase road throughput by driving at closer distances \cite{Lioris2017PlatoonsRoads} and can reduce fuel consumption (and therefore CO$_2$ emissions) up to $20\%$ \cite{Liang2016Heavy-DutyEfficiency}.
For both cooperative and automated driving (CAD), ensuring safety for higher levels of automation requires 
 architectures that contain health monitoring and management, safety-channels and fallback functionalities  \cite{ISO26262-2018,KhabbazSaberi2015277}.
The safety mechanisms designed to mitigate potential safety-critical hazards should be able to transition and bring a vehicle to a \textit{safe state}, i.e. an operating mode without an unreasonable level of risk. 
In addition, vehicles operating in SAE level 4 or 5 should be able to autonomously reach a minimal risk condition in case of a performance-relevant system failure \cite{TaxonomyAD}.
This implies that the vehicle should, without the interference of a human driver, bring itself to a minimal risk or safe condition when a fault or failure within the vehicle occurs, such that the vehicle can no longer be operated in the absence of unreasonable risk (e.g., a brake or steering failure, in the absence of any redundant or other risk mitigating measures).

To address the concerns above, the authors of \cite{Luo2017AnAnalysis} proposed an architecture pattern with a safety channel suitable for automated driving applications and Automotive Safety Integrity Level (ASIL) D, which is the highest risk class. 
In this work, the safety channel is divided into a health channel and a limp home channel. However, it does not specify the functionalities and methods that should be used in these channels, as this would highly depend on the level of automation and the type of functionalities involved. Falling back to such a channel would be a logical consequence of being able to diagnose a fault, crossing a level of severity which disables the vehicle to operate in a nominal condition. This requires functionalities to diagnose the system and check for the presence of faults and their severity. 
The survey \cite{Gao2015AApproaches} provides an overview of methods that can be used to \emph{diagnose} a fault, which implies three steps: (i) detection, i.e.  determining whether there is a fault, (ii) isolation, i.e. the location of the fault; (iii) identification, i.e the type, shape and size of the fault.
Furthermore, \cite{Gao2015AApproaches} briefly discusses fault tolerant control (FTC) strategies, where the system performance is maintained in the presence of faults, yet no real connection is made between the diagnosis and what mitigation measures should be taken.
Similarly, \cite{Yang2020Fault-tolerantMethodologies} fault tolerant cooperative control is introduced, focusing on mitigation strategies. All current work focuses on single mild faults, that require a limp home mode or degraded functionality and so far, there is no end-to-end system including functional safety considerations for both diagnosis and mitigation of faults of different types.

Once a severe fault is diagnosed, it is of foremost importance to bring the system to a safe state (i.e., make use of a safe and efficient fallback strategy). 
To this end, \cite{Svensson2018SafeFormulation} focuses on trajectory planning in fallback scenarios by formulating the problem as an optimal control problem, without considering any faults.
The authors of~\cite{Xue2018AEnvironment} describe an adaptive model predictive control (MPC) algorithm to simultaneously avoid potential collisions with surrounding vehicles and handle the presence of a front perceptive sensor failure.
Yu and Luo propose \cite{Yu2019FallbackSystem} a fallback strategy to park on the road shoulder while having a loss of all redundant paths or GPS location.
However, they decouple the longitudinal and lateral control of the vehicle, which might hinder the safe vehicle movement towards the road shoulder, especially in high risk scenarios.
Furthermore, in all previous works, no failures are considered that influence the vehicle's handling.

\textit{Our contributions:} the literature on trajectory planning and control for high vehicle automation is rich. Yet, it lacks work on fallback strategies for these functionalities, which form an essential part of the safety mechanisms in functional safety. We sum up our contributions as follows.
\begin{enumerate}[label=(\roman*)]
    \item The first contribution of this paper is a fallback strategy which is proven, based on the architectural \textit{pattern} proposed in~\cite{Luo2017AnAnalysis}. We pick up the architectural \textit{pattern} from~\cite{Luo2017AnAnalysis} and design a functionality with a software architecture that fits in this proven pattern, through which we accommodate fault diagnosis as well as mitigation for automated driving applications. The architectural design aims to facilitate all required steps from nominal operation to transitioning the vehicle to a safe-state in case of severe failures. 
    \item The second contribution of this paper is focusing on a vehicle affected by a failure and model uncertainty, for which an MPC-based fail-safe mitigation algorithm is introduced with coupled longitudinal and lateral dynamics. This algorithm, deployed inside a safety channel, ensures the safe operation of the vehicle in case of a failure, by bringing it to the emergency lane. 
\end{enumerate}

\begin{figure}[t]
    \centering
    \includegraphics[width = \linewidth]{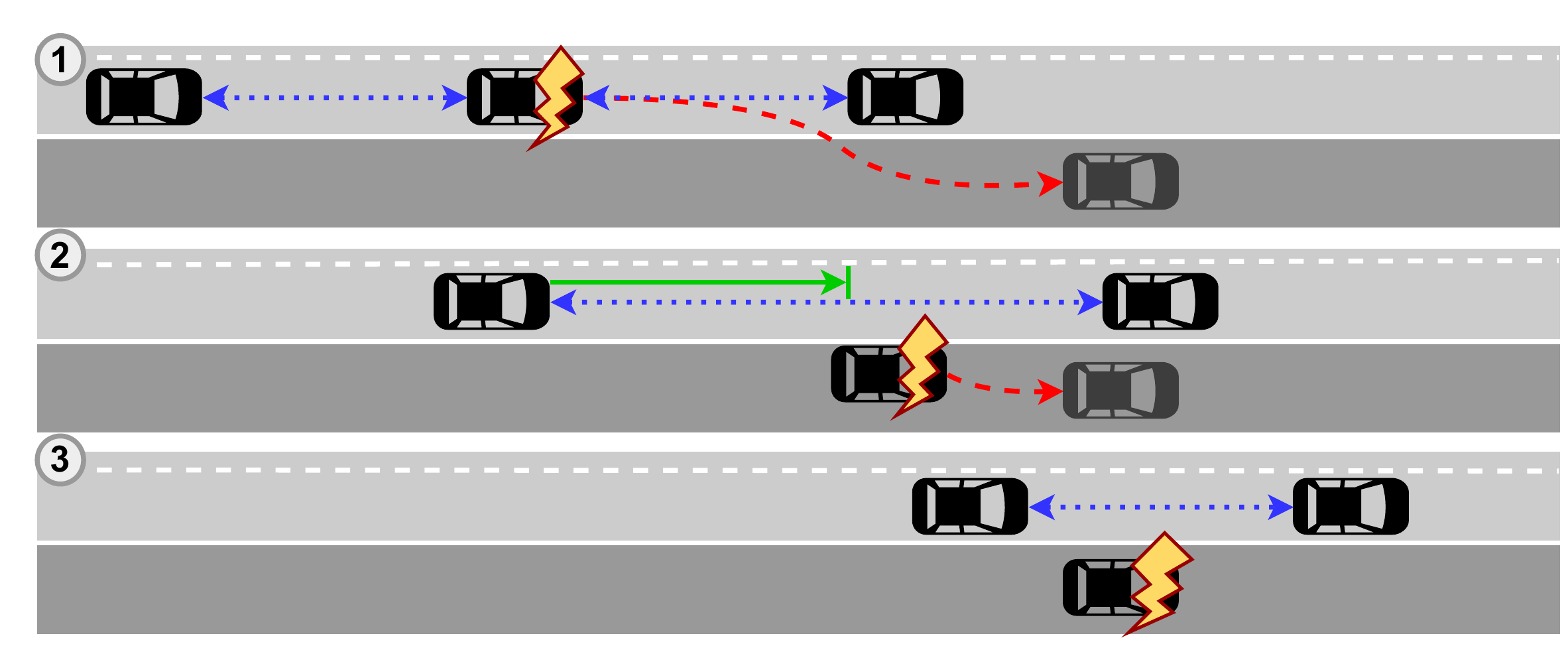}
    \caption{Scenario description where a severe fault occurs in a string of automated vehicles.}
    \label{fig:scenario}
\end{figure}

Fig. \ref{fig:scenario} shows the example scenario considered, where in a string of automated vehicles, running in nominal conditions, one detects a fault and needs to automatically park itself on the road shoulder. 
Within this scenario, two mitigation strategies are investigated for the faulty vehicle to showcase the influence on the remainder of the string of vehicles:
(i) The vehicle will brake inside the current lane, starting from the point that it receives the instruction to park on the road shoulder, and 
(ii) The vehicle will brake outside of the current lane, starting from the point that it has left the active lane.

This paper is organised as follows. Section II introduces the main components of the proposed architecture which enables nominal and fallback functionalities for an automated vehicle. Section III introduces the fail-safe mitigation algorithm and Section IV presents the simulation results for various scenarios and fault severity levels. Finally, conclusions and recommendations are described in Section V.

\section{FUNCTIONAL SAFETY ARCHITECTURE} \label{chap:Prop_Arch}
\begin{figure}[b]
    \centering
    \includegraphics[width = 0.8 \linewidth]{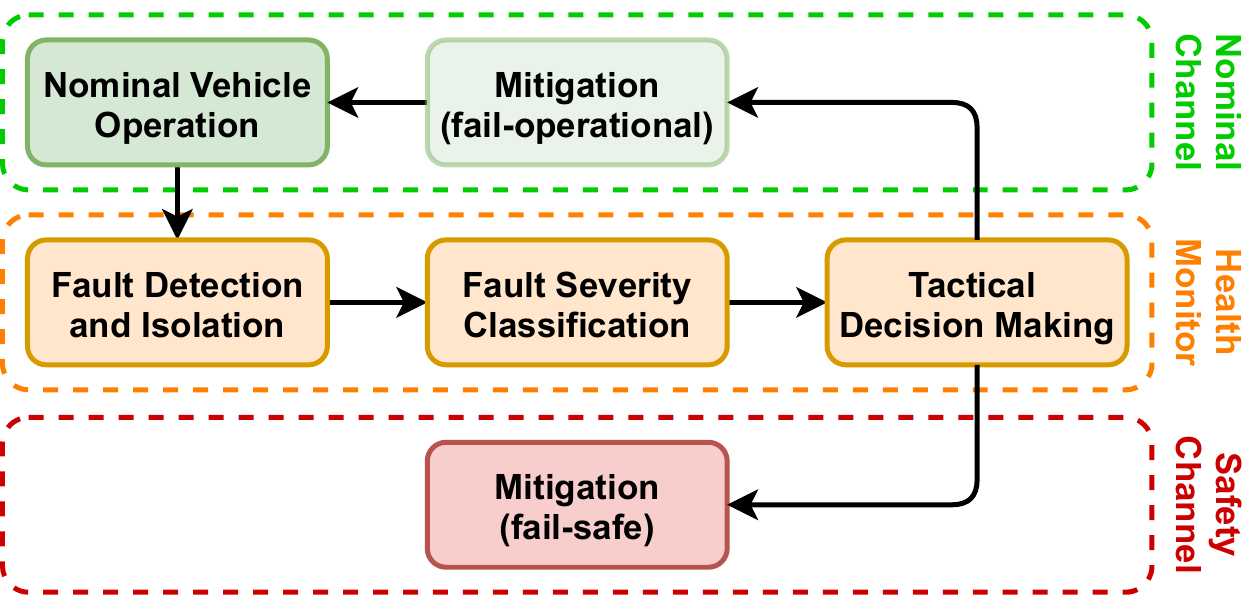}
    \caption{Architecture approach for nominal and safety fallback functionality of an automated vehicle, including failures.}
    \label{fig:Block_diag_Architecture}
\end{figure}

To ensure safe and comfortable operations, an automated vehicle architecture consists of three parts, namely, a nominal channel, a health monitor and a safety channel~\cite{Luo2017AnAnalysis}. We propose here the architecture shown in  Fig.~\ref{fig:Block_diag_Architecture}, which is an actual applied architecture based on the architectural pattern proposed in~\cite{Luo2017AnAnalysis}. 
Herein, the nominal channel performs all the nominal vehicle operation, i.e., all automated tasks which could function in the absence of unreasonable risk. The health monitor continuously monitors data coming from the vehicle to check whether this is operating in a healthy state, and the safety channel accommodates fail-safe mitigation to bring the vehicle to a safe-state when needed. The design of the actual module is not in the scope of this paper, however, earlier results show the feasibility of designing a suitable fault estimator~\cite{9756936}, which, given appropriate thresholds, can serve as a suitable classification algorithm.

\subsection{Nominal Vehicle Operation}
Nominal vehicle operation refers to the operation of the vehicle under normal circumstances, that is, in the absence of anomalies, faults or failures (AFFs) which could impose unreasonable risk to the vehicle and passengers. Moreover, in nominal operation, the vehicle is assumed to be driven in its Operational Design Domain. In these conditions, the system can make use of all its functionalities and ensure safe vehicle control. \looseness=-1

\subsection{Fault Detection and Isolation}
To assess AFFs, first, their presence and location should be known. 
This is done by respectively the detection and the isolation, where the detection solely focuses on the presence of an AFF. 
Subsequently, the isolation then determines the location of the AFF.
Finally, the identification determines its type, shape, and size, using advanced observer techniques such as Proportional (Multiple-) Integral observers, adaptive observers, sliding mode observers or descriptor observers \cite{Gao2015AApproaches}.

\subsection{Fault Severity Classification}
The risk of the diagnosed fault can be classified using ASIL levels and safety channel hardware, to determine if it is safe for the vehicle to continue driving. 
If it is not safe for the vehicle, the module determines whether the vehicle can continue with degraded functionality or whether it should go to a safe state.

\subsection{Tactical Decision Making}
In literature, this module is implemented both on a single- and multi-vehicle level, if the vehicle has vehicle-to-vehicle communication and can drive in cooperative modes (e.g., platooning \cite{konstantinopoulou2019specifications}). Tactical decision making is usually a needed nominal functionality that also contains a health monitoring and management component.
In CAD, by using this module, the integrity of a string of vehicles can be maintained while, for example, one of the vehicles abruptly leaves the string of vehicles. In the context of the scenario described in Fig. \ref{fig:scenario}, a benefit of this module is that the behaviour of the Lead Vehicle (LV) can be influenced such that the Trailing Vehicle (TV) can reconnect to the LV while optimizing certain parameters (e.g., fuel consumption).

\subsection{Mitigation}
Reducing the effect of an AFF is referred to as \textit{mitigation}.
Anomalies can lead to faults and consequently to failures, which are undesirable and potentially unsafe. In the context of a failure, i.e., a termination of an intended behaviour of an element or an item due to a fault manifestation \cite{ISO26262-2018}, handling this failure means controlling the system in its presence. Depending on the outcome of the Fault Severity Classification (FSC) module, the strategy for the mitigation is chosen to be fail-operational or fail-safe.
\paragraph{Fail-operational}
When the FSC module determines that the vehicle can safely continue operation, possibly with reduced functionality (also referred to as degraded or limp functionality), fail-operational mitigation is performed.
Such mitigation is most commonly performed by FTC if the AFF concerns an actuator or process \cite{Yang2020Fault-tolerantMethodologies}.
As exemplified in \cite{Khalili2018}, FTC converts the system to be less or not at all dependent on the faulty component, using the information acquired in the health monitor.
\paragraph{Fail-safe}
In case the FSC module determines that the vehicle is in a non-healthy state and cannot guarantee safe operation, fail-safe mitigation is performed by initiating a fallback manoeuvre to bring the vehicle to a safe-state. 
Similar to fail-operational mitigation, the information acquired in the health monitor is used.

\section{Fail-Safe Mitigation Algorithm} \label{chap:PoRS}

\begin{figure}[b]
    \centering
    \includegraphics[width = 0.9\linewidth]{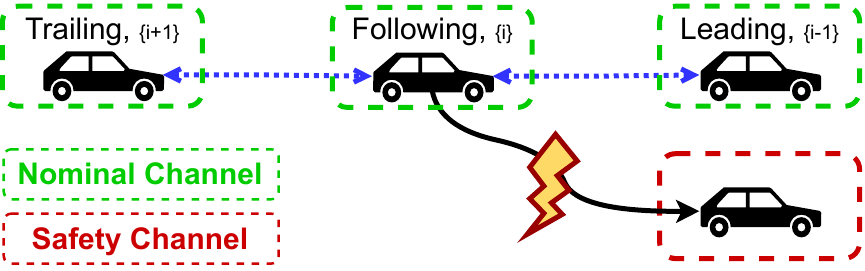}
    \caption{Multiple ACC-driven vehicles of which one encounters a severe fault and  needs to reach a safe state.}
    \label{fig:Control_distribution}
\end{figure}

To describe the fail-safe mitigation algorithm proposed in this paper, we start from the scenario described in Fig.~\ref{fig:scenario}.
Herein, three vehicles are assumed to drive automatically on the road (with functionalities such as adaptive cruise control and lane keep assist active, i.e., the \textit{nominal} functionality). As shown in Fig. \ref{fig:Control_distribution}, once a severe fault occurs, the faulty vehicle needs to transition to a safe state with the help of its safety channel.
The ACC-based longitudinal controller is ensuring a constant time-gap inter vehicle distance, with the following error dynamics
\begin{align}
    e_{tg} = h_{dg} - \frac{d_{x,i-1} - d_{x,i}}{v_{x,i}},
\end{align}
where $e_{tg}$ represents the time gap error between the two vehicles, $h_{dg}$ indicates the desired time gap between the two vehicles, $d_{x,i-1} - d_{x,i}$ is the distance between the preceding vehicle and the ego vehicle, and $v_{x,i}$ is the ego vehicle velocity.  \looseness=-1 

This error is controlled by a Proportional Derivative (PD) controller~\cite{PloegACC} with the control law formulated in the Laplace domain as follows:
\begin{equation}
    u_{PD} = e_{tg} (k_p  + k_d s),
\end{equation}
where $u_{PD}$ is the control output, $k_p$ the proportional gain, and $k_d$ the derivative gain of the control law.

To ensure safe handling of the faulty vehicle, both longitudinal and lateral control is immediately taken over by the safety channel after AFF diagnosis.
Without loss of generality, we assume here the faults are already detected and classified and focus on the Tactical Decision Making (TDM) and Fail-Safe Mitigation (FSM) modules from Fig. \ref{fig:Block_diag_Architecture}.
The implemented TDM is explained in Section \ref{sec:ImpTDM} and the controller used in FSM is explained in Section \ref{sec:Control}.

\subsection{Implemented Tactical Decision Making} \label{sec:ImpTDM}
\begin{figure}[t]
    \centering
    \includegraphics[width = 0.85\linewidth]{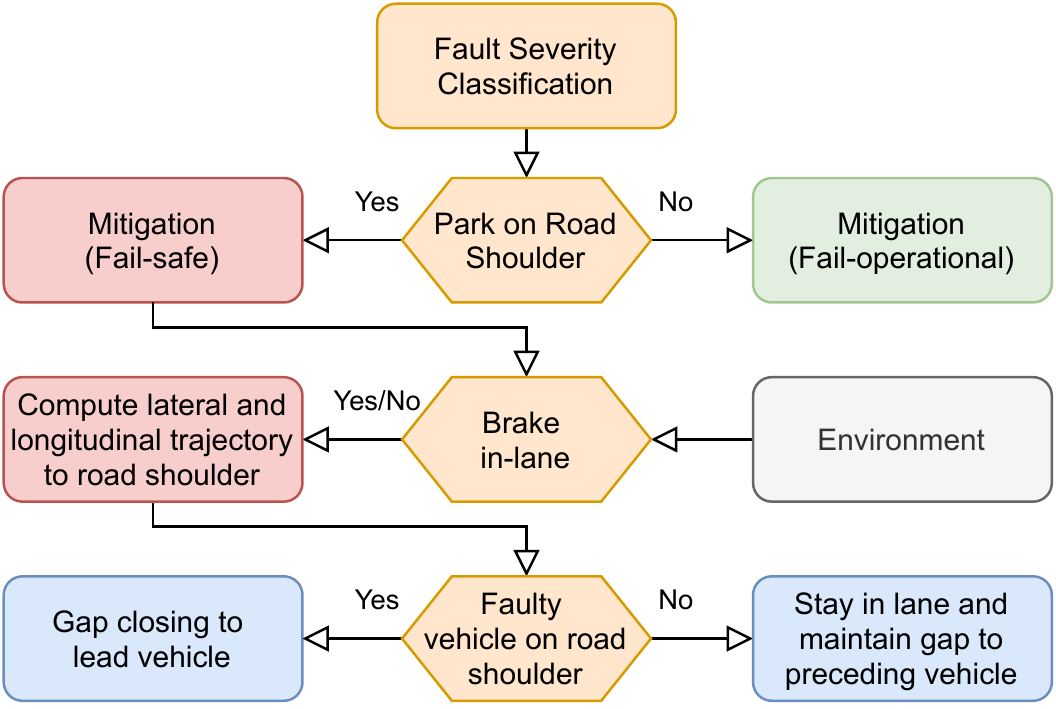}
    \caption{Flow chart of the implemented vehicle Tactical Decision Making (hexagons) and its effects on the trailing vehicle (blue blocks).}
    \label{fig:tactical_layer}
\end{figure}

Fig. \ref{fig:tactical_layer} shows the implemented TDM module, in which the FSC module gives a message to the TDM module when the failure is classified and the vehicle should be parked on the road shoulder.
The environmental module gives input that determines if the vehicle should brake in, or out-of-lane, e.g. if the road shoulder is long enough to brake out-of-lane, otherwise brake in-lane is required.
Eventually, the TDM module sends a message to the TV when it should close the gap back to the LV.
\label{sec:TacticalDecisionMaking_FSM}

\subsection{Functional Safety Mitigation Controller} \label{sec:Control}
MPC is often used to generate optimal control commands for the vehicle, \cite{MPC_explain, mpc_V2V, Ploeg:2022}, by taking into account the vehicle dynamics and its limitations over a predefined time window, known as prediction horizon $N$.
A Nonlinear MPC (NMPC) performs the high-level control in the safety channel and is required because of the combined longitudinal and lateral dynamics, based on the continuous-time equations of the linear single-track dynamic bicycle model~\cite{Schmeitz2017}:
\begin{align}
&\begin{aligned} \label{eq:vy_dot} 
&\dot{v}_y(t) = - \frac{C_{\alpha f} + C_{\alpha r}f_2}{m v_x(t)}v_y(t) ~+ \\
& \Bigg(\frac{l_r C_{\alpha r}f_2 - l_f C_{\alpha f}}{m v_x(t)} - v_x(t) \Bigg)r(t) + \frac{C_{\alpha f}}{m} \delta(t)f_1 
 \end{aligned} \\
&\begin{aligned}\label{eq:r_dot}
\dot{r}(t) &= \frac{l_r C_{\alpha r}f_2 - l_f C_{\alpha f}}{I_z v_x(t)}v_y(t) ~-\\
& \frac{l^2_f C_{\alpha f} + l_r^2 C_{\alpha r}f_2}{I_z v_x(t)}r(t) + \frac{l_f C_{\alpha f}}{I_z} \delta (t)f_1
\end{aligned}
\end{align}
Where $C_{\alpha f}$ and $C_{\alpha r}$ are the front and rear cornering stiffness, respectively, $m$ is the vehicle mass, $l_f$ and $l_r$ are the length from the front and rear axles to the center of gravity, respectively, $I_z$ is the vehicle's moment of inertia and finally, $v_x, \, v_y, \, r$ represent the longitudinal velocity, lateral velocity and yaw rate, respectively. Finally, $f_1$ and $f_2$ represent signals, acting on the steering wheel angle $\delta(t)$ and the rear corner stiffness $C_{\alpha r}$, respectively. These signals represent a power steering failure, $f_1$, and model uncertainty, $f_2$, potentially introduced by a fault. Throughout this paper, we assume that the signals $f_1,\:f_2$ appear as constants and are measurable by a fault diagnosis algorithm.
Within the constraints imposed by the linear bicycle model, the tyre dynamics are also linear.

The proposed MPC design requires a discrete-time update model, thus Equations \eqref{eq:vy_dot} and \eqref{eq:r_dot} are discretized using the forward Euler method, to form the nonlinear state-update equations (Equation \eqref{eq:mpc_prediction_model}), from the state vector $x(k)$:
\begin{align}
    x(k) = [a_x(k) \, v_x(k) \, v_y(k) \, d_y(k) \, r(k) \ \theta(k)]^\textrm{T},
\label{eq:State_vector}
\end{align}
where $a_x, \, d_y$ and $\theta$ are the longitudinal acceleration, and lateral position with respect to the center of the current lane and heading angle, respectively.
\begin{subequations}
\label{eq:mpc_prediction_model}
\begin{align}
    &a_x(k+1) = s_{dt} a_x(k) + G_{dt} a_{x,c}(k) \label{eq:ax_up} \\
    &v_x(k+1) = v_x(k) + a_x(k)\Delta t \label{eq:vx_up} \\
    &v_y(k+1) = v_y(k) + \Delta v_y(k) \Delta t \label{eq:vy_up} \\
    &\begin{aligned} d_y(k+1) &= d_y(k) ~+ \\ &\big(v_y(k)\cos\big(\theta(k)\big) + v_x(k)\sin\big(\theta(k)\big)\big)\Delta t \end{aligned} \label{eq:dy_up}\\
    &r(k+1) = r(k) + \Delta r(k) \Delta t \label{eq:r_up} \\
    &\theta (k+1) = \theta(k) + r(k) \Delta t \label{eq:theta_up}, 
\end{align}
\end{subequations}
where $s_{dt}$ and $G_{dt}$ are respectively the discrete-time pole and gain of the first-order transfer function representing the longitudinal dynamics, $a_{x,c}$ is the intended longitudinal acceleration, $\delta$ is the front wheel angle, $\Delta v_y$ and $\Delta r$ are the increments in $v_y$ and $r$, $\Delta t$ denotes the sampling time step and the indicator $k$ denotes the discrete time step. 

Model~\eqref{eq:mpc_prediction_model} can be rewritten in a more compact notation as:
\begin{equation}
    x(k+1) = g(x(k),u(k),f),
\end{equation}
Where $u(k) := \left [ a_{x,c},\, \delta\right ]^{\textrm{T}}$ is the control vector and $f:=\left[f_1,\:f_2\right]$ represent the constant values of the determined failure. The NMPC is formulated as
\begin{subequations}
\label{eq:NMPC}
\begin{align} 
    \min_u &\sum_{k=1}^{N} J\left(x(k),u(k),z(k)\right)\\
    \text{s.t. \hspace{0.5cm}} &x(k+1) = g\left(x(k),u(k),f\right) \label{eq:nmpc_dynamics}\\
    &x_{\min} \leq x(k) \leq x_{\max} \label{eq:xineq}\\
    &u_{\min}  \leq u(k) \leq u_{\max} \label{eq:uineq}\\
    &\Delta u_{\min} \leq \frac{u(k+1)-u(k)}{\Delta t} \leq \Delta u_{\max} \label{eq:duineq}\\
    &a_{y,\min} \leq a_y(k) \leq a_{y,\max}\label{eq:lateral_acc} \\
    &x(0) = x_{\textrm{init}}\\
    &\forall k \in \{0,\dots,N\},
\end{align}
\end{subequations}
where $z(k)$ contains the reference from the trajectory generation and $J$ represents the multi-objective cost function:
\begin{equation} \label{eq:cost_function}
\begin{aligned}
    J\Big(x(k),u(k),z(k)\Big) = w_{v_x}\big(z_{v_x}(k) - v_{x}(k)\big)^2 &+ \\ 
    w_{d_y}\big(z_{d_y}(k) - d_{y}(k)\big)^2 + 
    w_{\theta}\big(z_{\theta}(k) - \theta(k)\big)^2 &+ \\
    w_{a_x} \big(a_{x,c}(k)\big)^2 + w_{\delta} \big( \delta(k)\big)^2,
\end{aligned}
\end{equation}
where $w_{(\dots)}$ are the respective weights.
Constraint~\eqref{eq:nmpc_dynamics} indicates the dynamic coupling and constraints \eqref{eq:xineq}, \eqref{eq:uineq} and \eqref{eq:duineq} indicate comfort and model limitations.
Within which $\delta$, $\dot{\delta}$ and $\dot{a}_{x,c}$ are based on the physical capabilities of the vehicle and limits of the dynamic bicycle model. The constraints on $a_{x}, \, a_{x,c}, \, a_y$ and $v_x$ are based on the maximum allowed ACC braking, according to ISO 15622, comfort and highway speed limit respectively.
The lateral acceleration $a_y$ in~\eqref{eq:lateral_acc} is calculated by the following steady-state relation (imposed as a comfort constraint):
\begin{equation} \begin{aligned} \label{eq:ineqcon_ay}
    a_y &= -\frac{C_{\alpha f} + C_{\alpha r}}{m v_x}v_y \!+\! 
\frac{l_r C_{\alpha r} - l_f C_{\alpha f}}{m v_x} r \!+\! \frac{C_{\alpha f}}{m} \delta
\end{aligned} \end{equation}
Note, that the problem is assumed feasible in the scope of this work. However, through the use of an additional slack variable in the objective and carefully chosen constraints, one can enforce feasibility by sacrificing certain vehicle-dynamic constraints.
\section{SIMULATION RESULTS} \label{chap:Sim}

For this simulation study, a string of vehicles is considered as depicted in Fig.  \ref{fig:Control_distribution}, where all vehicles are modelled using the parameters given in Table \ref{tab:vehic_par}.a.
These parameters correspond to a lab passenger vehicle available at TNO~\footnote{https://www.tno.nl/en/focus-areas/traffic-transport/expertise-groups/research-on-integrated-vehicle-safety/}, used for research on cooperative and automated driving technologies.
The constraint values used in the NMPC model are given in Table \ref{tab:constraints}.b.

\begin{table}[t]
    \caption{Simulation parameters}
    \label{tab:vehic_par}
    \subfloat[Vehicle parameters]{
        \centering
        \vspace{0.14in}
        \begin{tabular}{c|c|c}
             Parameter & value & unit  \\ \hline \hline
            $C_{\alpha f}$ & 120 & $kN/rad$\\ 
            $C_{\alpha r}$ & 220 & $kN/rad$\\
            $l_f$ & 1.33 & $m$\\ 
            $l_r$ & 1.47 & $m$\\ 
            $m$ & 1845 & $kg$ \\ 
            $I_z$ & 3580 & $kg \cdot m^2$
        \end{tabular}
            \vspace{0.18in}
        
    }\label{tab:constraints}
    \subfloat[Constraints parameters]{
        \centering
        
        \begin{tabular}{c|c|c}
            Variable & Constraint & unit \\
            & (min / max) & \\ \hline \hline
            $\mid \delta \mid$ & 0.0873 & $rad$ \\
            $\mid \dot{\delta} \mid$ & 0.0818 & $rad/s$ \\
            $a_{x}$ & -3.5 / 1.5 & $m/s^2$ \\
            $a_{x,c}$ & -3.5 / 1.5 & $m/s^2$ \\
            $\dot{a}_{x,c}$ & -14 / 6 & $m/s^3$ \\
            $v_x$ & 1.26 / 33 & $m/s$ \\
            $\mid a_y \mid$ & 2 & $m/s^2$
        \end{tabular}
    }
\end{table}
The trajectory that the faulty vehicle follows during the fail-safe mitigation is split into lateral and longitudinal movement, to best accommodate both our mitigation strategies.
The lateral trajectory is generated by a $\text{5}^{\text{th}}$ order polynomial, taken between current and goal waypoints with appropriate heading angles, following \cite{Yu2019FallbackSystem}.
The current waypoint is the middle of the active lane and the goal waypoint is the middle of the road shoulder, assuming a straight road.
For the longitudinal trajectory, only goal velocities are given, such that the controller determines the optimal control outputs within the given constraints, considering all relevant dynamics.
Alternatively, as part of our future work, a local motion planner can also be incorporated into our architecture to adapt the trajectory online to avoid collisions with upcoming traffic (e.g., \cite{safeVRU}).

The vehicle model that is used as a plant, to test the controller, is based around the continuous time counterparts in \eqref{eq:mpc_prediction_model}.  \looseness=-1

\subsection{Controller settings}
The tuning parameters for the Proportional Derivative (PD) controllers performing the longitudinal control for the ACC string of vehicles and the NMPC controller that performs the fallback manoeuvre are given in Tables \ref{tab:PD_control} and \ref{tab:NLMPC_control}, respectively. \looseness=-1

Table \ref{tab:PD_control} shows the tuning parameters of each vehicle, where the LV is tuned differently compared to the FV and TV, as it is operating in cruise control and tracking a reference velocity instead of a time-gap to the preceding vehicle. \looseness=-1
\begin{table}[b]
    \centering
    \caption{Settings of PD controllers of each vehicle}
   \begin{tabular}{c|c|c}
         Vehicle & $k_p$ & $k_d$ \\  \hline \hline
         Leading & 5 & 0.3 \\
         Following~/~Trailing & -150 & -2.5 
    \end{tabular}
     \label{tab:PD_control}
\end{table}

From the dynamic bicycle model in~\eqref{eq:vy_dot} and~\eqref{eq:r_dot} it can be derived that, as the velocity decreases towards zero, the eigenvalues of the linear differential equations grow towards $-\infty$. This phenomenon is numerically impossible to capture in the Forward Euler approximation used in this paper, as it would require the sampling time to be reduced to $0$. Following this line of reasoning, to prevent numerical instability of the internal prediction model, a sampling time of $0.01\,s$ and a $v_{x,\min}$ of $1.26\,m/s$ is selected.
The selection of the prediction horizon $N$, control horizon $S$ and the NMPC weights $w_{(\dots)}$ in the cost function are manually chosen with a trade-off between computational effort and tracking performance, aiming for low computational effort with minimal loss in tracking performance. \looseness=-1
\begin{table}[t]
    \centering
    \caption{Settings of NMPC used for the fallback manoeuvre}
    \begin{tabular}{c|c|c|c|c|c|c|c} 
         Variable & N & S & $w_{v_x}$ & $w_{d_y}$ & $w_{\theta}$ & $w_{a_x}$ & $w_{\delta}$  \\ \hline \hline
        Value & 30 & 30 & 10 & 100 & 1 & 0.5 & 1
    \end{tabular}
    \label{tab:NLMPC_control}
\end{table}

\subsection{Failure scenarios and braking strategies considered} \label{sec:FailScen}

We present six simulation results: (i) two simulations compare braking in-lane and braking out-of-lane during the fallback manoeuvre, (ii) two simulations investigating the robustness of the controller by implementing realistic failures and uncertainties in the vehicle model and (iii) two simulations investigating the behaviour of the controller if it is reconfigured, following the architecture proposed in Section \ref{chap:Prop_Arch}, adjusting relevant formulas and bounds.

The following failure and uncertainty are considered for (ii) and (iii):
\begin{enumerate}
    \item $f_1$: Power steering failure
The steering output of the controller is decreased by $50 \%$ before it feeds through to the vehicle model, thus $f_1 = 0.5$. 
    \item $f_2$: Model uncertainty in the rear cornering stiffness
The rear cornering stiffness $C_{\alpha r}$ of the vehicle is decreased by $50 \%$, thus $f_2 = 0.5$.
\end{enumerate}

\subsection{Results}
To show the performance of the proposed method in bringing the vehicle to a safe state, two time moments are important, $t_a$, when the parking manoeuvre is initiated, and $t_b$, when the faulty vehicle has left the initial driving lane. 

\subsection{Braking in-lane versus braking out-of-lane}
Table \ref{tab:Result_values} highlights the trade-off between the two mitigation strategies based on stop time and distance versus the re-connection time of the remaining vehicles on the road (TV to the LV).
The stop time is calculated as the time between $t_a$ and the time that the error on the goal velocity is less or equal to $0.01\,m/s$ and the error on the lateral position is less or equal to $0.001\,m$.
The travelled distance between these two instances is the stopping distance.
Re-connection time is calculated as the time between the instance that $e_{tg}$ is larger than $0.4\,s$ and the instance that the $e_{tg}$ stays below $0.01\,s$.
\begin{table}[t]
    \caption{Comparison between braking strategies while going to the road shoulder and the effect on upcoming traffic.}
    \begin{tabular}{c|c|c|c|c}
    Braking   & Stop      & Stop          & Trailer      & Timegap \\
    mitigation     & time [s]  & distance      & gap-closing         & error\\
    strategy & & [m] &time [s]  & at $t_b$ [s]\\\hline \hline 
    In-lane & ~8.208 & 117.534 & 13.880 & 1.650 \\ 
    Out-of-lane & 10.838 & 190.610 & ~7.634 & 1.004
    \end{tabular}
    \label{tab:Result_values}
\end{table}
Stopping time and distance are largely influenced by $a_{x,\min}$ as the lateral movement consumes less time compared to the longitudinal movement.
Furthermore, the duration of the lateral movement has a major impact on the difference in closing time due to the distance and velocity difference it creates between both strategies.

\begin{figure}[b]
    \centering%
    \includegraphics[clip, trim=0cm 0cm 0cm 0cm]{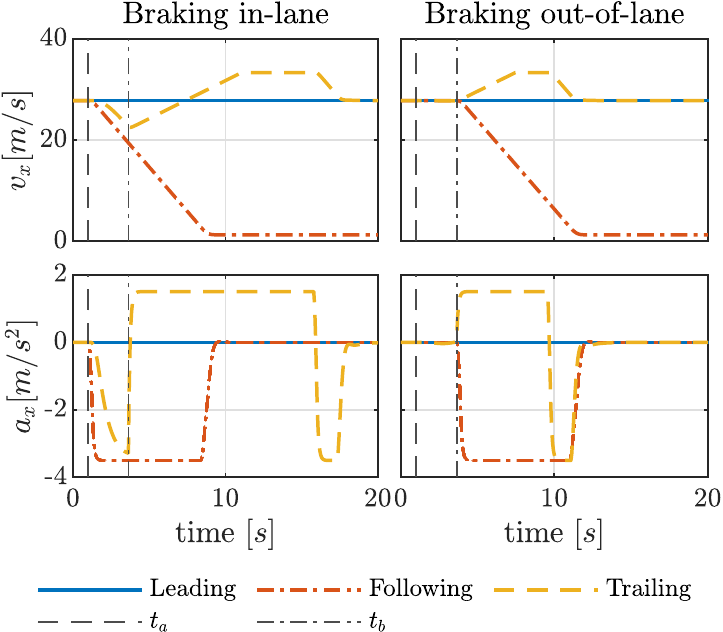}
    \caption{Longitudinal velocities $v_x$ and accelerations $a_x$ during the lane changing strategies for all vehicles.}
    \label{fig:axvxBILBOL}
\end{figure}

As expected, the timegap error at $t_b$ shows that braking in-lane (BIL) results in a higher time-gap than braking out-of-lane (BOL) and therefore a longer closing time for BIL compared to BOL.
This is underpinned by the velocity difference between the TV and LV at $t_b$ and the acceleration length in Fig. \ref{fig:axvxBILBOL}.
\begin{figure}[t]
    \centering
    \includegraphics[clip, trim=0cm 0cm 0cm 0cm]{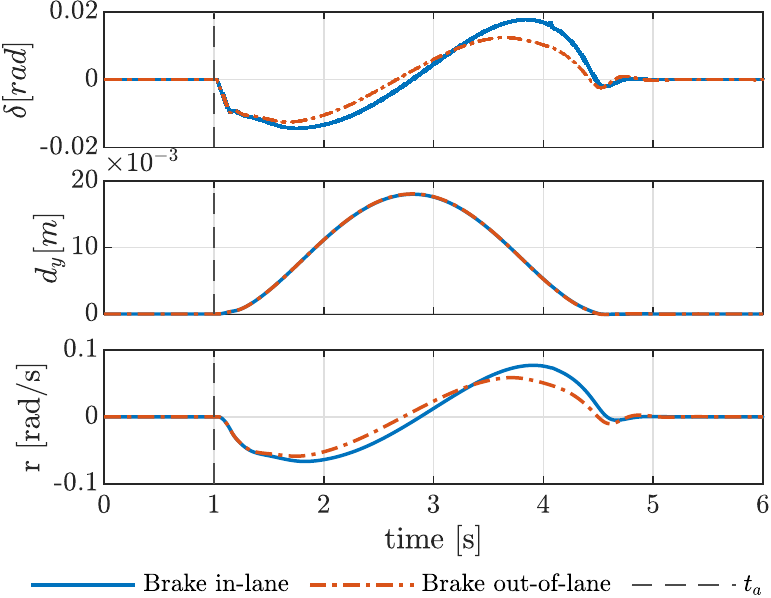}
    \caption{Comparison between both mitigation strategies on steering output $\delta$, lateral position $d_y$ and yaw rate $r$}
    \label{fig:BILBOLvehicle}
\end{figure}
The lateral deviation in both strategies is equal, following Fig. \ref{fig:BILBOLvehicle}, however, the steering outputs show different behaviour in both strategies.
This, helped by the decreased longitudinal velocity because of braking, translates into an increased yaw rate in the vehicle dynamics for BIL compared to BOL.

As BIL results in higher lateral loads on the vehicle dynamics, this strategy is used in further experiments and as a baseline comparison.
Figs. \ref{fig:Failures} and \ref{fig:Failures_comp} show the error difference between the input/states of the baseline (BIL without failure) and the input/states of the subsequent failure.

\subsection{Robustness of the controller}
\begin{figure}[t]
    \centering
    \includegraphics[clip, trim=0cm 0cm 0cm 0cm]{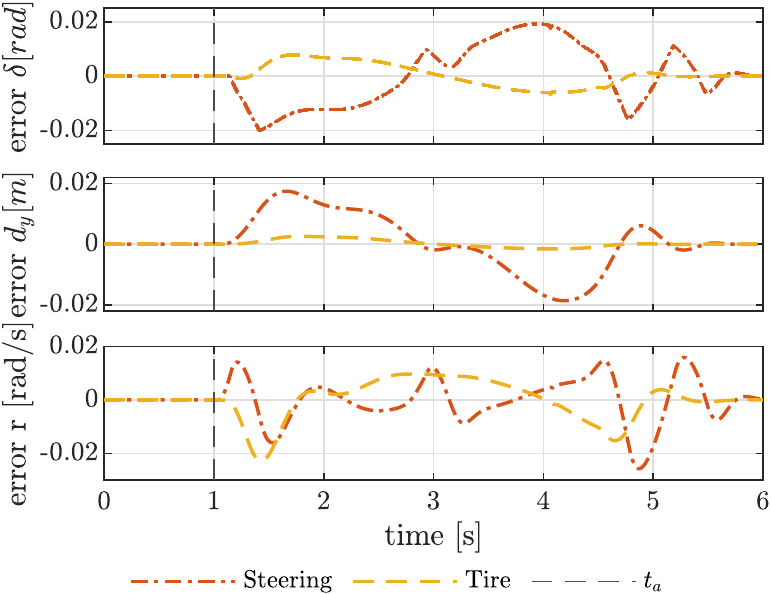}
    \caption{Error plots of the controller with a steering failure and uncertainty in the rear cornering stiffness on steering output $\delta$, lateral position $d_y$ and yaw rate $r$ compared to the baseline.}
    \label{fig:Failures}
\end{figure}

Following the results in Fig. \ref{fig:Failures}, the steering failure causes the controller to output higher steering inputs for the vehicle. Next to that, steering is less smooth and shows more abrupt changes in direction, caused by reaching the limit of the steering rate $\dot{\delta}$.
This is also clearly visible in the yaw rate $r$, showing its influence on the lateral vehicle dynamics.

Furthermore, the results show that the model uncertainty decreases the maximum controller setpoint and makes the initially understeered vehicle show oversteered behaviour.
The latter translates into the vehicle turning more compared to the baseline with the same steering input.
This effect is also observed in the lateral position error, where the vehicle initially steers too much, and thus deviates further from the path.

\subsection{Reconfiguration of the controller}
\begin{figure}[t]
    \centering
    \includegraphics[clip, trim=0cm 0cm 0cm 0cm]{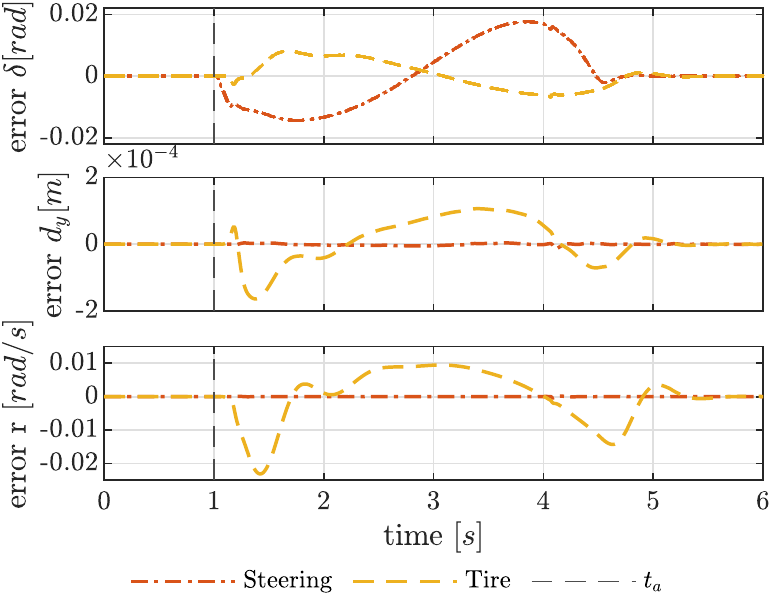}
    \caption{Error plots of the reconfigured controller with a steering failure and uncertainty in the rear cornering stiffness on steering output $\delta$, lateral position $d_y$ and yaw rate $r$ compared to the baseline.}
    \label{fig:Failures_comp}
\end{figure}

The reconfigured controller uses the information on the failures (Section \ref{sec:FailScen}) to update the internal NMPC model (Equation \eqref{eq:mpc_prediction_model}).
For the steering failure, this means that $\delta(k)$ is transformed into $0.5 \delta (k)$.
Also, the bounds on $\delta$ and $\dot{\delta}$ are increased by a factor $\frac{1}{0.5}$.
In the case of the model uncertainty in the rear cornering stiffness, $C_{\alpha r}$ is changed to half of its original value.
Fig. \ref{fig:Failures_comp} shows the results, in which the lateral control action is smoother for the steering failure but similar for the model uncertainty, compared to the non-reconfigured controller in Fig. \ref{fig:Failures}.
The magnitude of the steering output is comparable to the failure and the model uncertainty in relation to the non-reconfigured simulations.

The steering output and yaw rate of the model uncertainty can be compared with its non-reconfigured result, however, the maximum lateral deviation error is decreased by $92 \%$ for the reconfigured controller.
For the steering failure, the error on lateral deviation is decreased to under $0.013 \, mm$, a decrease up to $33 \%$, and the yaw rate error to a maximum of $0.00037 \, rad/s$, effectively eliminating the effect of the failure on tracking performance.

When evaluating all figures and  Table \ref{tab:Result_values}, it is clear that the controller is capable of handling a power steering failure or model uncertainty in the rear cornering stiffness.
Especially when re-configuring the NMPC model, the performance is comparable to that of the system without failure.
Furthermore, as BOL results in a lower gap-closing time for the TV, thus disrupts the surrounding vehicles less than BIL, and results in lower dynamic loads thus higher comfort, it is recommended to use this mitigation strategy if the environment of the vehicle allows this.  \looseness=-1

\section{CONCLUSIONS} \label{chap:Conclusion}
The contributions of this research focuses on introducing a functional safety architecture that can handle multiple types of faults, the strategy and the fail-safe mitigation algorithm to park the vehicle on the road shoulder in case of severe failures. Such an architecture is essential to enable higher levels of automation and prove the functional safety of a system when a failure occurs.
  
Our fail-safe mitigation strategy (tactical decision making and motion control) relies on a finite state machine and a tailored MPC formulation, controlling the lateral and longitudinal movement of the vehicle simultaneously.
The results, shown for a severe failure (i.e. power steering failure) and model uncertainty in the rear cornering stiffness, highlight the trade-offs for different lane changing strategies for the faulty vehicle, i.e. braking in- and out-of-lane, and for the other vehicle in upcoming traffic.
Furthermore, results also show that if the controller can have failure-awareness it can adapt and performance can be improved. 

In future work we plan to validate our proposed architecture and fail-safe mitigation algorithm also through experiments, to verify it using more scenarios and by incorporating the other needed components (such as fault diagnosis and severity classification).
Furthermore, we aim to perform a stability analysis on the proposed NMPC controller and look further into the consequences on the remainder of the platoon (e.g. on string stability and time headways). Other work includes real-time implementation and experimental validation.

\section{Acknowledgements}
This work is supported by the EU Horizon 2020 R{\&}D program under grant agreement No. 861570, project SAFE-UP (proactive SAFEty systems and tools for a constantly UPgrading road environment) and from the Dutch Science Foundation NWO-TTW, within the Veni project HARMONIA (nr. 18165).

\errorcontextlines=99

\end{document}